\documentclass[a4paper,10pt]{article}
\usepackage[scaled]{helvet}              
\usepackage[T1]{fontenc}                 
\usepackage{geometry}  
\usepackage{setspace}                 
\usepackage{exscale}                  
\usepackage{a4wide}                   
\usepackage{framed}                   
\usepackage{subscript}                
\usepackage{authblk}                  
\usepackage{appendix}                 
\usepackage{colordvi}                 
\usepackage{enumitem}                 
\usepackage[american]{babel}
\usepackage{fancyhdr}                   
\usepackage[colorlinks=true,linkcolor=black,citecolor=black,filecolor=black,urlcolor=black,hypertexnames=false]{hyperref} 
\usepackage{aliascnt}      
\usepackage{cleveref}      
\usepackage[babel,german=guillemets]{csquotes}
\usepackage[backend=bibtex,style=numeric-comp,useprefix=true,hyperref=true,firstinits=true]{biblatex}   
\bibliography{../../../Import/References/references.bib}
\usepackage[fleqn]{amsmath} 
\usepackage{mathtools}      
\usepackage{amssymb}        
\usepackage{latexsym}       
\usepackage{dsfont}         
\numberwithin{equation}{section}      
\input{../../../Import/Environments/thm-environments.tex}
\input{../../../Import/Macros/mathmacros-global.tex}

\usepackage{xargs}   
\newcommand{\spaceTC}{ \Hk{1}{} }                                                         
\newcommand{\spaceSC}{ \fspace{L^2}{I}{; \Hk{1}{} } \cap \fspace{H^1}{I}{; \Hk{1}{}^* } } 
\newcommand{\spaceTF}{ \Hkdiv[][0]{} }                                                    
\newcommand{\spaceSF}{ \fspace{L^\infty}{I}{; \Hkdiv[][f]{} } }                           
\newcommand{\spaceTP}{ \Lp{2}{} }                                                         
\newcommand{\spaceSP}{ \fspace{L^\infty}{I}{; \Lp{2}{}/\setR } }                          
\newcommand{\spaceSEF}{ \fspace{L^\infty}{I}{; \Hkdiv[][\sigma]{} }  }                    
\newcommand*{\tx}{(t,x)}    
\newcommand*{\timeEnd}{T_0} 

\newcommand*{\test}{\ensuremath{ \vphi }} 
\newcommand*{\Test}{\ensuremath{ \vecv }} 

\newcommand*{\boltz}{\ensuremath{ k_b }}                                                 
\newcommand*{\temp}{\ensuremath{ T }}                                                    
\newcommand*{\Dl}[1][]{\ensuremath{\mathcal{D}_{{#1}} }}                                 
\newcommand*{\zl}[1][l]{\ensuremath{ z_{{#1}} }}                                         
\newcommand*{\permitELscalar}{\ensuremath{ \epsilon }}                                   
\newcommand*{\permitEL}{\ensuremath{ \mathcal{E} }}                                      
\newcommand*{\permeabH}{\ensuremath{ \mathcal{K} }}                                      
\newcommand*{\coeffEL}[1][l]{\ensuremath{e\zl[{#1}](\permitELscalar\boltz\temp)^{-1} }}  
\newcommand*{\coeffForceEL}{\ensuremath{\mu^{-1}\permitEL^{-1}}}

\newcommand{\cl}[1][l]{\ensuremath{ c_{{#1}} }}                  
\newcommand{\Cl}{\ensuremath{ \vecc }}                           
\newcommand{\clstart}[1][l]{\ensuremath{ c_{0,{#1}} }}           
\newcommand{\Rl}[1][l]{\ensuremath{ R_{{#1}} }}                                 

\newcommand{\fieldF}{\ensuremath{ \vecu }} 
\newcommand{\press}{p}                     

\newcommand{\potEL}{\ensuremath{ \Phi }}                                    
\newcommand{\fieldEL}{\ensuremath{ \vecE }}                                 
\newcommand*{\chargeEL}[1][]{ \ensuremath{ \rho_{f{#1}} } }                 
\newcommand*{\chargeELlong}[1][l]{ \sum_l\zl\cl }                           
\newcommand*{\chargeELlongTwo}[1][l]{ \zl[1]\cl[1]-\abs{\zl[2]}\cl[2] }     
\newcommand*{\forceEL}[1][]{\ensuremath{ \chargeEL[{#1}] \fieldEL_{{#1}} }} 
\newcommand{\fixOp}{\ensuremath{\mathcal{F}}}
\newcommand{\fixSet}{\ensuremath{K}}

\newcommand{\clold}[1][l]{\ensuremath{\bar{c}_{{#1}} }}                           
\newcommand{\Clold}{\ensuremath{ \bar{\vecc} }}                                   
\newcommand*{\chargeELold}[1][]{ \ensuremath{ \bar{\rho}_{f{#1}} } }              
\newcommand*{\chargeELlongold}[1][l]{ \sum_l\zl\clold }                           
\newcommand*{\chargeELlongTwoold}[1][l]{ \zl[1]\clold[1]-\abs{\zl[2]}\clold[2] }  
\newcommand*{\forceELold}[1][]{\ensuremath{ \chargeELold[{#1}] \fieldEL_{{#1}} }} 
\newcommand{\lyapunov}{\Lambda}

\input{../../../Import/Environments/environments.tex}
\newcounter{countAssumption} 
\newcounter{countNotation}   
\newcommand{\labelInt}{I.a}
\newcommand{\labelNot}{(N1)--(N4)}
\begin{document}
%
%
\title{Global existence of weak solutions of a model for electrolyte solutions -- Part 2: Multicomponent case}
\author[1]{Matthias Herz} 
\author[1]{Peter Knabner}
\renewcommand\Affilfont{\itshape\small}
\affil[1]{Department of Mathematics, University of Erlangen-N\"urnberg, Cauerstr. 11, D-91058 Erlangen, Germany}
\date{\today}
\maketitle 
%
%
\begin{abstract}
This paper analytically investigates the \dpnp. This system is a mathematical model for electrolyte solutions. In this paper, we consider electrolyte solutions, 
which consist of a neutral fluid and multiple suspended charged chemical species with arbitrary valencies. We prove global existence and uniqueness of weak solutions in two space dimensions.
\\[2.0mm]
\end{abstract}
%
%
\section{Introduction}\label{sec:Introduction}
From the introduction of Part 1 of this work, we know that many complicated phenomena in hydrodynamics and biology can be modeled as electrolyte solutions, as these models simultaneously capture the following 
three ubiquitous processes: (i) The transport of the charged particles, (ii) the hydrodynamic fluid flow, (iii) the electrostatics.  
\medskip
\par
\dpnp{s} are classical models for electrolyte solutions in case of laminar flow in porous media. A detailed derivation of these systems are given in 
\cite{Allaire10, Elimelech-book, Juengel-book, Masliyah-book, Russel-book, Samohyl, Schuss, Mielke10, Probstein-book, RayMunteanKnabner}. 
In Part~1 of this work, we mentioned, that among many others, mathematical models for electrolyte solutions have been investigated analytically in 
\cite{BT2011, Burger12, Castellanos-book, GasserJungel97, Glitzky2004, Gajewski, HyongEtAll, Markovich-book, Mielke10, Roubicek2005-1, Schuss, Wolfram, Schmuck11, herz_ex_dnpp}. 
\par
In particular, existence of solutions for electrolyte solutions, which contain multiple charged solutes, were proven amongst others in \cite{BT2011, Burger12, Glitzky2004, Roubicek2005-1, Roubicek2006}. 
More precisely, the authors of \cite{BT2011} proved local in time existence for a fluid-particle system. The existence results of \cite{Roubicek2005-1, Roubicek2006} were proven under the additional assumption 
of a volume-additivity constraint and by including an additional reaction force term in the transport equations. These additional assumptions allow to bypass in \cite{Roubicek2005-1, Roubicek2006} 
the below mentioned serious difficulties. \cite{Burger12} dealt with a stationary system and in \cite{Glitzky2004}, existence in two dimensions was established for electrolyte solution at rest, 
which may involve nonlocal constraints in the equation for the electrostatic potential.
\medskip
\par    
In case of multiple charged chemical species with arbitrary valencies, the main difficulty is to prove a~priori estimates for the chemical species, which are valid on arbitrary large time intervals 
and which are independently of the electric field. 
We briefly sketch this difficulty by firstly considering the case of two oppositely charged species~$\cl[1]$ (positively charged) and $\cl[2]$ (negatively charged) 
with valencies~$\zl[1]>0>\zl[2]$. In Part~1 of this work, we proposed to use in the proof of a~priori estimates for $\cl$ the weighted test functions~$\test=\abs{\zl}\cl$. 
Together with Gauss's law, i.e. $\grad\cdot\fieldEL=\zl[1]\cl[1]-\abs{\zl[2]}\cl[2]$, we thereby obtained the following pointwise sign condition 
for the sum of the \enquote{electric drift integrals}, which describe the electrophoretic motion of the~$\cl$
\begin{align*}
 &-2\sum_l \scp[]{\zl\cl\fieldEL}{\grad\brac{\abs{\zl}\cl}~}_{\Lp{2}{}} = \scp[]{\grad\cdot\fieldEL}{\brac{\zl[1]\cl[1]}^2-\brac{\abs{\zl[2]}\cl[2]}^2}_{\Lp{2}{}} \\
 &= \scp[]{\zl[1]\cl[1]-\abs{\zl[2]}\cl[2]}{\brac{\zl[1]\cl[1]}^2-\brac{\abs{\zl[2]}\cl[2]}^2}_{\Lp{2}{}}  \geq 0, ~~\text{as } \sqbrac{a - b}\sqbrac{a^2 - b^2} \geq 0 ~\text{ for } a,b\geq0~.
\end{align*}
Due to this pointwise sign condition, we omitted the sum of the \enquote{electric drift integrals} in the proof of a~priori estimates for the chemical species~$\cl$.
However, in the multicomponent case, i.e., in case of  of $L\in\setN$ charged solutes~$\cl$ with arbitrary valencies~$\zl$, the corresponding sum of these integrals reads 
with Gauss's law~$\grad\cdot\vecE=\sum_l\zl\cl$ as
\begin{align*}
 &-2\sum_l \scp{\Sign(\zl) \abs{\zl}\cl \vec{E}}{\grad(\abs{\zl}\cl)~}_{\Lp{2}{}} = \scp{\grad\cdot\vecE}{\sum_l \Sign(\zl))(\abs{\zl}\cl)^2~ }_{\Lp{2}{}} \\
 &= \scp{ \sum_l \Sign(\zl)\abs{\zl}\cl }{\sum_l  \Sign(\zl)(\abs{\zl}\cl)^2}_{\Lp{2}{}}~~\Longleftrightarrow:~ I.1 = I.2 = I.3.
\end{align*}
It is easy to show that $I.3$ does not satisfy a pointwise sign condition%
\footnote{Consider, e.g., three chemical species with valencies~$\zl[1]=\zl[2]=1$, $\zl[3]=-1$. Assume that the values at a given point in space and time are $\cl[1]=\cl[2]=1$, $\cl[3]=\sqrt{3}$. This gives $(1+1-\sqrt{3})(1+1-3)=-(2-\sqrt{3})<0$.}. 
Thus, we now have to estimate one of the integrals~$I.1$--$I.3$ by suitable norms of the integrands. This can easily be done with standard arguments. More precisely, we obtain with \Holder, \GagNirenberg, and \Young 
the estimates (space dimension~$n=2,3$) 
\begin{align*}
 (a)~~& I.3 \leq C \sum_l \norm{c_l}{\Lp{3}{}}^3  \leq \delta \sum_l \norm{\grad\cl}{\Lp{2}{}}^2 + \delta^{-1}C \sum_l \norm{\cl}{\Lp{2}{}}^6~,  \\     
 (b)~~& I.2 \leq \delta \sum_l \norm{\grad\cl}{\Lp{2}{}}^2  + \delta^{-1}C \norm{\grad\cdot\vecE}{\Lp{2}{}}^{4/(4-n)} \sum_l \norm{\cl}{\Lp{2}{}}^2, \\
 (c)~~& I.1 \leq \delta \sum_l \norm{\grad\cl}{\Lp{2}{}}^2  + \delta^{-1}C \norm{\vecE}{\Lp{6}{}}^4 \sum_l \norm{c_l}{\Lp{2}{}}^2~.
\end{align*}
Provided we have an a~priori bound for the electric field or its divergence, which is independent of $\cl$, we are done by choosing one of the estimates~$(b)$ or $(c)$. 
However, in \pnp{s} the equations for the electric field and its divergence depend on $\cl$, which leads to the a~priori bound $\norm{\vec{E}}{\Lp{6}{}} + \norm{\grad\cdot\vec{E}}{\Lp{2}{}} \leq C \sum_l \norm{c_l}{\Lp{2}{}}$. 
Substituting this bound into $(b)$ or $(c)$ shows that these estimates lead to $(a)$. Finally, the proof of a~priori estimates for $\cl$ in combination with $(a)$ results in
\begin{align*}
 \derr[t]\sum_l \norm{c_l}{\Lp{2}{}}^2 + \sum_l \norm{\grad c_l}{\Lp{2}{}}^2 \leq C \sum_l \norm{c_l}{\Lp{2}{}}^2 + C \sum_l \norm{c_l}{\Lp{2}{}}^6~.
\end{align*}
Hence, we can not apply the standard version of \Grownwall, which we need to obtain a uniform a~priori bound from the preceding equation. In fact, we just have a nonlinear version due to 
Willett-Wong, cf. \cite{willett_wong}, which is in our case applicable just on a very small time interval, cf. \cite[Theorem~4.9]{Bainov-book}. An iterative application of this nonlinear \Grownwall\ does not lead 
to an extension to arbitrary large time intervals$[0,\timeEnd]$, as the constants may blow up in finite time.  
\medskip
\par
In this paper, we propose to substitute Gauss's law~$\grad\cdot\fieldEL=\sum_l\zl\cl=:\chargeEL$ into the above estimate~$(b)$. Thereby, estimate~$(b)$ reads as  
\begin{align*}
 I.2 &\leq \delta \sum_l \norm{\cl}{\Hk{1}{}}^2 + \delta^{-1} C\norm{\chargeEL}{\Lp{2}{}}^{4/(4-n)} \sum_l \norm{\cl}{\Lp{2}{}}^2~.
\end{align*}
It now remains to show a uniform bound for $\chargeEL$ in $\Lp[I]{4/(4-n)}{;\Lp{2}{}}$, which is the crucial task in the current chapter. We successfully establish such a uniform bound by means of entropy estimates, 
which lead to bounds for certain Lyapunov functions and which have been previously used in a different context, e.g., in \cite{KrautleExPaper}. However, these techniques lead to a uniform bound 
for $\chargeEL$ in $\Lp[I]{2}{;\Lp{2}{}}$. This is the reason why the presented existence result is restricted to two space dimensions.    
\medskip
\par
The contribution of this paper is to show in two space dimensions the global existence, uniqueness, and boundedness of weak solutions for electrolyte solutions that consist of a neutral fluid and 
multiple charged solutes with arbitrary valencies. Due to the entropy estimates we can avoid further restrictions, such as the electroneutralitiy constraint, cf. \cite{Allaire10}, or 
the volume-additivity constraint, cf. \cite{Roubicek2005-1}. Thus, the presented results apply to two dimensional models of general electrolyte solutions, which are captured by the \dpnp.
\medskip
\par
The rest of this paper is organized as follows: In \cref{sec:Model}, we present the model equations and we prove that solutions are unique. In \cref{sec:fpi-operator}, 
we introduce the fixed point method, and in \cref{subsec:aprioriBounds}, we prove the crucial a~priori estimates. Finally, in \cref{subsec:fpi-point}, we show the global existence.
Note, that we omit several proofs in \cref{subsec:aprioriBounds} and \cref{subsec:fpi-point}, as we wrote Part~1 of this work such that it allows us to copy the coinciding parts of the proof. 
%
%
\section{Model Equations}\label{sec:Model}
Henceforth, we use the notation \labelNot, which we already introduced in Part~1 of this work. Furthermore, in Part~1 of this work, we pointed out that \dpnp{s} consist of the following three 
coupled conservation laws:\\[1.0mm]
\begin{subequations}
%
\textbf{Law 1 -- Gauss's law} 
\vspace{-2mm}
\begin{align}
              \permitEL^{-1}\fieldEL   &=  -\grad\potEL                &&  \text{in } \OmegaT,   \label{eq:Model-1a} \\
           \grad\cdot\fieldEL          &=   \chargeEL       + \rho_b   &&  \text{in } \OmegaT,   \label{eq:Model-1b} \\
                           \chargeEL   &=  \theta\chargeELlong         &&  \text{in } \OmegaT,   \label{eq:Model-1c} \\[-3.0mm]
             \fieldEL\cdot\vecnu       &=  \sigma                      &&  \text{on } \GammaT,   \label{eq:Model-1d}
\end{align}
%
\textbf{Law 2 -- Darcy's law}
\vspace{-2mm}
\begin{align}
    \permeabH^{-1} \fieldF    &= \mu^{-1}\brac{-\grad\press + \permitEL^{-1}\forceEL}  &&  \text{in } \OmegaT,  \label{eq:Model-1e}  \\
	  \grad\cdot\fieldF   &=  0                                                    &&  \text{in } \OmegaT,  \label{eq:Model-1f}  \\
	  \fieldF\cdot\vecnu  &=  f                                                    &&  \text{on } \GammaT.  \label{eq:Model-1g} 
\end{align} 
%
\textbf{Law 3 -- Nernst-Planck equations} 
\vspace{-2mm}
\begin{flalign}
  \theta\dert\cl + \grad\cdot\brac{\Dl\grad\cl + {\cl}[\fieldF + \coeffEL\fieldEL] } &= \theta\Rl(\Cl)   &&  \text{in } \OmegaT,  & \label{eq:Model-1h} \\[0.2cm]
                  \brac{\Dl\grad\cl + {\cl}[\fieldF + \coeffEL\fieldEL] }\cdot\vecnu &= 0                &&  \text{on } \GammaT,  & \label{eq:Model-1i} \\[0.2cm]
                                                                              \cl(0) &= \clstart         &&  \text{on } \Omega.   & \label{eq:Model-1j}
\end{flalign}
\end{subequations}
\begin{remark}\label{remark:fluxBoundary}
\textnormal{
The homogeneous boundary flux in equation~\eqref{eq:Model-1i}, can be equivalently expressed with equations \eqref{eq:Model-1d} and \eqref{eq:Model-1g} by
\begin{align*}
 \brac{\Dl\grad\cl}\cdot\vecnu &= - {\cl}f - (\permitELscalar\boltz\temp)^{-1}e\zl{\cl}\sigma  \qquad \text{on }~ \GammaT~.
\end{align*}
Thus, the homogeneous boundary flux condition is equivalent to a homogeneous Robin boundary condition for the diffusion part.
}\hfill$\square$
\end{remark}
%
%
\subsection{Weak formulation of the model}\label{subsec:weakFormulation}
Firstly, we introduce the required notation for the analytical investigations. 
%
\begin{enumerate}[label=({N}\arabic*), ref=({N}\arabic*), itemsep=0.0mm, start=5]
 \item \textbf{Spaces: } For $k>0$, $p\in [1,\infty]$, we denote the Lebesgue spaces for scalar-valued and vector-valued functions by $\Lp{p}{}$ 
       and the respective Sobolev spaces by $\Wkp{k}{p}{}$, cf.~\cite{Adams2-book}.  
       Furthermore, we set $\Hk{k}{}:=\Wkp{k}{2}{}$ and we refer for the definition of the Bochner spaces~$\fspace{L^p}{I}{;V}$, $\fspace{H^k}{I}{;V}$ 
       over a Banach space~$V$ to~\cite{Roubicek-book}. The $\Hkdiv[][f]{k}$-spaces are defined, e.g., in~\cite{brezzi-book} by \\ 
       $\Hkdiv[][f]{k}:=\cbrac{\Test\in \Hk{k}{}:~\nabla\cdot\Test\in \Hk{k}{}, \Test\cdot\nu = f \text{ on } \Gamma}$.%
       \label{Not:spaces} 
 \item \textbf{Products: } We denote by $\scp{\cdot}{\cdot}_H$ the inner product on a Hilbert space~$H$ and by $\dualp{\cdot}{\cdot}_{V^\ast\times V}$, the dual pairing between 
       a Banach space~$V$ and its dual space~$V^\ast$. On $\setR^n$, we just write $\vecv\cdot\vecu:=\scp{\vecv}{\vecu}_{\setR^n}$ and on $\Lp{2}{;\setR^d}$, 
       we just denote $\scp{\cdot}{\cdot}_{\Omega}:=\scp{\cdot}{\cdot}_{\Lp{2}{;\setR^d}}$. In particular the dual pairing between $\Hk{1}{}$ and its dual~$\Hk{1}{}^\ast$, 
       we abbreviate by $\dualp{\cdot}{\cdot}_{1,\Omega}:=\dualp{\cdot}{\cdot}_{\Hk{1}{}^\ast\times\Hk{1}{}}$.%
       \label{Not:Prod}
\setcounter{countNotation}{\value{enumi}}
\end{enumerate}
Secondly, we henceforth impose the following assumptions:
%
\begin{enumerate}[label=({A}\arabic*), ref=({A}\arabic*), itemsep=0.0mm,  start=1]
 \item \textbf{Geometry: } Let $n=2$ and $\Omega\subset\setR^n$ be a bounded Lipschitz domain, i.e. $\Gamma\in C^{0,1}$.%
       \label{Assump:Geom}
 \item \textbf{Initial data: } The initial data~$\clstart$ are non negative and bounded, i.e., \\
       $0\leq\clstart(x)\leq M_0$ for a.e. $x\in\Omega$ for some $M_0 \in\setR_+$.%
       \label{Assump:InitData}
 \item \textbf{Ellipticity: } The diffusivity tensor~$\Dl$ and the permeability tensor~$\permeabH$ satisfy \\ 
       $\Dl\xi\cdot\xi>\alpha_D\abs{\xi}^2$ and $\permeabH^{-1}\xi\cdot\xi>\alpha_K\abs{\xi}^2$ for all $\xi\in\setR^n$, \\ 
       $\Dl\xi\cdot\eta<C_D\abs{\xi}\abs{\eta}$ and $\permeabH^{-1}\xi\cdot\eta<C_K\abs{\xi}\abs{\eta}$ for all $\xi,\eta\in\setR^n$.%
       \label{Assump:Ellip}
 \item \textbf{Coefficients: } The porosity~$\theta$, the dynamic viscosity~$\mu$, and the electric permittivity~$\epsilon$ are positive constants.%
       \label{Assump:Coeff}
 \item \textbf{Reaction rates: } The reaction rate functions $\Rl:\setR^L\rightarrow\setR$ are global Lipschitz continuous functions, i.e., $\Rl\in \Ck[\setR^L]{0,1}{}$ 
       with Lipschitz constant~$C_{\Rl}$. Furthermore, we assume $\Rl(\vecnull)=0$ and $\Rl(v_1,\ldots,v_l,\ldots,v_L) \geq 0$ for all $\vecv\in\setR^L$ 
       with $v_l\leq 0$. This means, in case a chemical species vanishes, it can only be produced.%
       \label{Assump:Reaction} 
 \item \textbf{Boundary data: } We assume $\sigma\in\Lp[\OmegaT]{\infty}{}$ and $f\in\Lp[\OmegaT]{\infty}{}$. 
       Furthermore, we suppose that functions $\vecf,\vecsigma\in\Lp[I]{\infty}{;\Wkp{1}{\infty}{}}$ with $\vecsigma\cdot\vecnu=\sigma$ and $\vecf\cdot\vecnu=f$ exist.%
       \label{Assump:BoundData}
 \item \textbf{Background charge density: } We assume $\rho_b\in \Lp[\OmegaT]{\infty}{}$.%
       \label{Assump:back-charge}
\setcounter{countAssumption}{\value{enumi}}
\end{enumerate}
Thirdly, we define the weak formulation of the \dpnp.       
%
\begin{definition}[Weak solution]\label{def:weaksolution}
The vector $\brac{\fieldEL,\potEL,\fieldF,\press,\Cl}\in \setR^{2+2n+L}$ is a weak solution of the \dpnp\ \eqref{eq:Model-1a}--\eqref{eq:Model-1j},  if and only if
\begin{subequations}
\begin{enumerate}[label=(\roman*), ref=(\roman*), itemsep=-1.5mm]
 \item $(\fieldEL,\potEL)\in \spaceSEF\times\spaceSP$ solve for all $(\Test,\test) \in \spaceTF\times\spaceTP$ 
      \begin{align}
        \scp{\permitEL^{-1}\fieldEL}{\Test}_\Omega &= \scp{\potEL}{\grad\cdot\Test}_\Omega  \label{eq:gaussWeak}\\
        \scp{\grad\cdot\fieldEL}{\test}_{\Omega} &= \scp{\rho_b + \chargeEL}{\test}_\Omega \label{eq:gaussWeakDIV}. 
      \end{align}
 \item $(\fieldF,\press)\in \spaceSF\times\spaceSP$ solve for all $(\Test,\test) \in \spaceTF\times\spaceTP$
       \begin{align}
           \scp{\permeabH^{-1}\fieldF}{\Test}_\Omega &= \scp{\mu^{-1}\press}{\grad\cdot\Test}_\Omega + \scp{\coeffForceEL\forceEL}{\Test}_\Omega, \label{eq:darcyWeak}\\
               \scp{\grad\cdot\fieldF}{\test}_\Omega &= 0~.    \label{eq:darcyWeakDIV}
       \end{align}
 \item $\cl \in \fspace{L^\infty}{I}{;\Lp{2}{}} \cap \spaceSC \cap \Lp[\OmegaT]{\infty}{}$ solves for all $\test\in\spaceTC$ and for $l=1,...,L$
       \begin{flalign}\label{eq:transportWeak}
        &~~\dualp{\theta\dert\cl}{\test}_{1,\Omega} + \scp{\Dl\grad\cl}{\grad\test}_\Omega - \scp{ {\cl}[\fieldF + \coeffEL\fieldEL] }{\grad\test}_\Omega 
          = \scp{\theta\Rl(\Cl)}{\test}_\Omega~,\!\!\!\!&
       \end{flalign}
       and $\cl$ take its initial values in the sense that 
       \begin{align*}
        \lim\limits_{t\searrow0} \scp{\cl(t) - \clstart}{\test}_\Omega ~=~0 \qquad \text{for all } \test\in\Lp{2}{}~.\\[-12mm] \nonumber
      \end{align*}
      \hfill$\square$
\end{enumerate}
\end{subequations}
\end{definition}
\begin{remark}
\textnormal{
As already showed in Part~1 of this work, equations~\eqref{eq:transportWeak} are not well-defined without having $\cl \in\Lp[\OmegaT]{\infty}{}$. 
Thus $\cl\in\Lp[\OmegaT]{\infty}{}$ is mandatory for a well-defined weak formulation.%
}\hfill$\square$ 
\end{remark}

%
%
\subsection{Uniqueness}\label{sec:Uniqueness}
We note that we can skip some of the following proofs, as we wrote Part~1 of this work such that we can directly copy some proofs. 
%
\begin{thm}[Uniqueness]\label{thm:unique}
Let \ref{Assump:Geom}--\ref{Assump:back-charge} be valid and let $\brac{\fieldEL,\potEL,\fieldF,\press,\Cl}\in \setR^{2+2n+L}$ be a weak solution 
of \eqref{eq:Model-1a}--\eqref{eq:Model-1j} according to \cref{def:weaksolution}. Then, $\brac{\fieldEL,\potEL,\fieldF,\press,\Cl}$ is unique. 
\end{thm}
\begin{proof}
The proof is identical to the corresponding proof of Part~1 of this work. Furthermore, the proof holds even in three space dimensions.
\end{proof}
%
%
\section{Fixed Point Operator} \label{sec:fpi-operator}
We apply the same fixed point approach as already used in Part~1 of this work. For more details concerning this approach, we refer to Part~1 of this work.  
%
%
\begin{definition}[Fixed point operator]\label{def:FixOp}
Let $K\subset X$ be a subset of the Banach space $X$, which is given by $X:= \sqbrac{\fspace{L^\infty}{I}{;\Lp{2}{}} \cap \spaceSC \cap \Lp[\OmegaT]{\infty}{}}^L$. 
We introduce the fixed point operator~$\fixOp$ by
\begin{align*}
\fixOp:=\fixOp_3\circ\fixOp_2\circ\fixOp_1: K \subset X \rightarrow X~.
\end{align*}
Herein, the suboperator $\fixOp_1$ is defined by
\begin{subequations}
\begin{align}
 \fixOp_1: &\begin{cases} \fixSet \rightarrow X\times\spaceSEF\times\spaceSP ~=:Y \\
                           \Clold ~\mapsto (\Clold,\fieldEL,\potEL), ~\text{ with } (\fieldEL,\potEL) \text{ solving for all } (\Test,\test) \in \spaceTF\times\spaceTP
              \end{cases} \nonumber\\
            & \hspace{27.0mm} \scp{\permitEL^{-1}\fieldEL}{\Test}_\Omega = \scp{\potEL}{\grad\cdot\Test}_\Omega , \label{eq:gaussWeak-fpi}    \\
            & \hspace{27.0mm} \scp{\grad\cdot\fieldEL}{\test}_{\Omega} = \scp{\rho_b + \chargeELold}{\test}_\Omega , \label{eq:gaussWeakDIV-fpi} \\
            & \hspace{27.0mm} \text{with } \chargeEL \text{ defined ind } \eqref{eq:Model-1c}~. \nonumber 
\end{align}
Furthermore, the suboperator $\fixOp_2$ is defined by
\begin{align}
 \fixOp_2: &\begin{cases} \qquad Y ~~~~\rightarrow Y\times \spaceSF\times\spaceSP~=:Z \\
                          (\Clold,\fieldEL,\potEL) \mapsto (\Clold,\fieldEL,\potEL,\fieldF,\press), ~\text{ with } (\fieldF,\press) \text{ solving for all } (\Test,\test) \in \spaceTF\times\spaceTP
              \end{cases} \nonumber\\
            & \hspace{27.0mm} \scp{\permeabH^{-1}\fieldF}{\Test}_\Omega = \scp{\mu^{-1}\press}{\grad\cdot\Test}_\Omega + \scp{\coeffForceEL\forceELold}{\Test}_\Omega, \label{eq:darcyWeak-fpi}   \\
            & \hspace{27.0mm} \scp{\grad\cdot\fieldF}{\test}_\Omega = 0 .   \label{eq:darcyWeakDIV-fpi} 
\end{align}
Finally, the suboperator $\fixOp_3$ is defined by
\begin{align}
 \fixOp_3: &\begin{cases} \qquad~~~ Z \qquad~~\rightarrow X \\
                          (\Clold,\fieldEL,\potEL,\fieldF,\press) \mapsto \Cl=(\cl[1],\cl[2]), ~\text{ with } \cl \text{ solving for all } \test\in\spaceTC \text{ and } l=1,\ldots,L
              \end{cases} \nonumber\\
            & \hspace{27.0mm} \dualp{\theta\dert\cl}{\test}_{1,\Omega} + \scp{\Dl\grad\cl}{\grad\test}_\Omega - \scp{ {\cl}[\fieldF + \coeffEL\fieldEL] }{\grad\test}_\Omega \nonumber\\[2.0mm]
            & \hspace{27.0mm} = \scp{\theta\Rl(\Cl)}{\test}_\Omega, \label{eq:transportWeak-fpi} \\
            & \hspace{27.0mm} \text{and } \cl \text{ take its initial values in the sense that } \nonumber\\
            & \hspace{27.0mm} \lim\limits_{t\searrow0} \scp{\cl(t) - \clstart}{\test}_\Omega ~=~0 \qquad \text{for all } \test\in\Lp{2}{}~. \nonumber\\[-12mm] \nonumber
\end{align}
\end{subequations}
\hfill$\square$
\end{definition} 
%
%
\begin{lemma}[well-definedness]\label{lemma:wellDef-fixOp}
Let \ref{Assump:Geom}--\ref{Assump:back-charge} be valid. Then, the operator $\fixOp: \fixSet\subset X \rightarrow X~$ defined in \cref{def:FixOp}, is well-defined.
\end{lemma}
\begin{proof}
 The proof is identical to the corresponding proof of Part~1 of this work. Furthermore, the proof holds even in three space dimensions.
\end{proof}
%
%
\begin{lemma}[regularity for Gauss's law]\label{lemma:regularity-gausslaw}
Let \ref{Assump:Geom}--\ref{Assump:back-charge} be valid and let $(\fieldEL,\potEL,\fieldF,\press,\Cl)\in\setR^{2+2n+L}$ be a solution of \eqref{eq:Model-1a}--\eqref{eq:Model-1j} according to \cref{def:FixOp}. 
Then, for the partial solution~$(\fieldEL,\potEL)$, we have
\begin{align*}
 \potEL\in\Lp[I]{\infty}{;\Hk{2}{}/\setR} \qquad \text{ and } \qquad \fieldEL\in\Lp[I]{\infty}{;\Hk{1}{}}~.
\end{align*}
\end{lemma}
\begin{proof}
 The proof is identical to the corresponding proof of Part~1 of this work. Furthermore, the proof holds even in three space dimensions.
\end{proof}
%
%
\section{Existence for multicomponent electrolytes} \label{sec:Existence-multComp}
In this section, we prove that global weak solutions of the \dpnp\ exist.
%
%
\subsection{A priori Estimates}\label{subsec:aprioriBounds}
We now show a~priori bounds for the solution vector $\brac{\fieldEL,\potEL,\fieldF,\press,\Cl}\in \setR^{2+2n+L}$. First of all, we cite some preliminary results, which we use in the subsequent calculations.
%
%
\begin{lemma}[Boundary Interpolation]\label{lemma:interpol-boundary}
 Let $u\in\Hk{1}{}$ and suppose \ref{Assump:Geom}. Then, we have
 \begin{align*}
  \norm{u}{\Lp[\Gamma]{2}{}}^2 ~\leq~ \delta \norm{\grad u}{\Lp{2}{}}^2 + 2\delta^{-1} \norm{u}{\Lp{2}{}}^2 \qquad \text{ for all } \delta\in(0,1)~.  
 \end{align*}
\end{lemma}
\begin{proof}
 The stated inequality follows immediately from \cite[Theorem~7.58]{Adams1-book}, \cite[Lemma~7.16]{Adams1-book}, and \Young.
\end{proof}
%
%
\begin{lemma}[Approximation]\label{lemma:approxInTime}
Let $I\subset\setR$ be an interval and let $u\in\Lp[I]{2}{;\Hk{1}{}}\cap\Hk[I]{1}{;\Hk{1}{}^\ast}$. Then, there exists a sequence $(u_\veps)_{\veps>0}\subset\Ck[I]{1}{;\Hk{1}{}}$ such that we have
\begin{align*}
 \norm{u-u_\veps}{\Lp[I]{2}{;\Hk{1}{}}}+\norm{u-u_\veps}{\Hk[I]{1}{;\Hk{1}{}^\ast}}+\norm{u-u_\veps}{\Lp[I]{\infty}{;\Lp{2}{}}} ~\rightarrow~ 0 ~.
\end{align*}
\end{lemma}
\begin{proof}
See \cite[Lemma~7.2, Lemma~7.3]{Roubicek-book}. 
\end{proof}
%
%
\begin{lemma}[Lyapunov function]\label{lemma:lyapunov}
Define the Lyapunov function~$\lyapunov(x):\setR_+\rightarrow\setR$ by 
\begin{align*}
 \lyapunov(x):=x(\ln(x)-1)+e \qquad\Hence \lyapunov^\prime(x) = \ln(x) ~.
\end{align*}
Then, we have  
\begin{align*}
 \lyapunov(x)-x &\geq 0 \qquad\text{and}\qquad \lyapunov(x)\geq 0 ~.
\end{align*}
\end{lemma}
\begin{proof}
 The second inequality follows with $x\ge0$ from the first one. For the first inequality calculate the minimum of the function~$\lyapunov(x)-x$. This shows $\min_{x\geq0}(\lyapunov(x)-x)=0$.
\end{proof}
\begin{remark}
We call the function~$\lyapunov(x)$ a Lyapunov function, since for the heat equation~$\dert u -\Delta u =0$ in weak formulation, we easily obtain by formally testing with $\ln(u)$ the estimate
\begin{align*}
 \norm{\lyapunov(u(t))}{\Lp{1}{}} ~\leq~ \norm{\lyapunov(u(0))}{\Lp{1}{}}~.
\end{align*}
Hence, the function $t\mapsto \norm{\lyapunov(u(t))}{\Lp{1}{}}$ is non-increasing along every trajectory~$t\mapsto u(t)$. From the theory of ordinary differential equations, 
we know that functions with this property are called Lyapunov functions, cf. \cite{Arnold-book}. 
\par
In the theory of partial differential equations, the Lyapunov function~$\lyapunov(x)$ occurs naturally in estimates that are based on testing with the logarithm of the solution. 
These estimates are called entropy estimates, cf. \cite{CarilloJuengelEtAl2001}, as in many situations, we can physically interpret for a solution~$u$ the logarithm~$\ln(u)$ 
as entropy, cf. \cite{LifshitzLandau-book5}. For that reason, the estimates in \cref{lemma:entropyBound} are called entropy estimates.
\par
Lyapunov functions of the type~$\lyapunov(x)$ have been used in \cite{KrautleExPaper} in the context of reaction-diffusion systems with reaction rates according to mass-action law.  
\hfill$\square$
\end{remark}
Next, we repeat a short result from Part~1 of this work 
%
%
\begin{lemma}[Non negativity]\label{lemma:nonnegative-Lcomp}
Let \ref{Assump:Geom}--\ref{Assump:back-charge} be valid and let $\brac{\fieldEL,\potEL,\fieldF,\press,\Cl}\in \setR^{2+2n+L}$ 
be a weak solution of \eqref{eq:Model-1a}--\eqref{eq:Model-1j} according to \cref{def:FixOp}. Then, we have for $l\in\{1,\ldots,L\}$
\begin{align*}
\cl\tx \geq 0 \quad \text{ for a.e. } t\in[0,\timeEnd],~ \text{ a.e. } x\in\Omega~.
\end{align*}
\end{lemma}
\begin{proof}
  The proof is identical to the corresponding proof of Part~1 of this work. Furthermore, the proof holds even in three space dimensions.
\end{proof}
We now prove the entropy estimates. Note, that these estimates continue to hold even for three space dimensions.
%
%
\begin{lemma}[Entropy estimates]\label{lemma:entropyBound}
Let \ref{Assump:Geom}--\ref{Assump:back-charge} be valid and let $\brac{\fieldEL,\potEL,\fieldF,\press,\Cl}\in \setR^{2+2n+L}$ be a solution of the system~\eqref{eq:Model-1a}--\eqref{eq:Model-1j} 
according to \cref{def:FixOp}. Then, we have the estimates
\begin{align*}
                                               \norm{\chargeELold}{\Lp[I]{2}{;\Lp{2}{}}} &\leq C_L(\timeEnd)~,\\
  \sum_l \norm{\sqrt{\cl}}{\Lp[I]{\infty}{;\Lp{2}{}}}^2 + \norm{\grad\sqrt{\cl}}{\Lp[I]{2}{;\Lp{2}{}}} &\leq C_L(\timeEnd)~.
\end{align*}
\end{lemma}
\begin{proof}
Let $a>0$. We test equations~\eqref{eq:transportWeak-fpi} with $\test:=\ln(\cl+a)\in\spaceTC$, we sum over $l$ and we integrate in time over $[0,t_1]\in[0,\timeEnd]$. 
For ease of readability, we split the following proof into two cases.\\[2.0mm]
%
\underline{Case 1: $\cl=\clold$}~~
For the sum of the time integrals, we obtain with \cref{lemma:approxInTime} and \cref{lemma:lyapunov} 
\begin{align*}
 &  \sum_l\dualp{\dert\cl}{\ln(\cl+a)}_{\Lp[0,t_1]{2}{;\Hk{1}{}^\ast}\times\Lp[0,t_1]{2}{;\Hk{1}{}}} \\
 &= \sum_l\lim\limits_{\veps\rightarrow0} \Intdt[][t_1]{\Intdx{ \dert(\cl[l,\veps]+a)\ln(\cl[l,\veps]+a)~}} 
 ~= \sum_l\lim\limits_{\veps\rightarrow0} \Intdt[][t_1]{\Intdx{ \dert\lyapunov(\cl[l,\veps]+a) ~}} \\ 
 &= \sum_l \norm{\lyapunov(\cl(t_1)+a)}{\Lp{1}{}} - \sum_l\norm{\lyapunov(\clstart+a)}{\Lp{1}{}}~.
\end{align*}
For the sum of the diffusion integrals, we immediately get with \ref{Assump:Ellip} 
\begin{align*}
      \sum_l \scp{\Dl\grad\cl}{\frac{1}{\cl+a}\grad(\cl+a)}_{\Omega\times[0,t_1]} 
 \geq \alpha_D \sum_l \norm{\grad\sqrt{\cl+a}}{\Lp[\Omega\times{[0,t_1]}]{2}{}}^2~.
\end{align*}
The sum of the convection integrals, we transform with integration by parts and equation~\eqref{eq:darcyWeakDIV-fpi} to
\begin{align*}
 I_{co}:=&   -\sum_l\scp{\cl\fieldF}{\grad\ln(\cl+a)}_{\Omega\times{[0,t_1]}} \\
         &=  -\sum_l\scp{(\cl+a)\fieldF}{\grad\ln(\cl+a)}_{\Omega\times{[0,t_1]}} - a\sum_l\scp{\fieldF}{\grad\ln(\cl+a)}_{\Omega\times{[0,t_1]}} \\
         &=  -\sum_l\scp{f}{\cl+a}_{\Gamma\times{[0,t_1]}} - \sum_l\scp{f}{a\ln(\cl+a)}_{\Gamma\times{[0,t_1]}}
         ~=: A.1 + A.2~.
\end{align*}
Applying \Holder\, \cref{lemma:interpol-boundary}, and \cref{lemma:lyapunov} leads for the integral~$I.1$ 
to
\begin{align*}
 A.1 &\geq -\sum_l\norm{f}{\Lp[\GammaT]{\infty}{}}\norm{\cl+a}{\Lp[{\Gamma\times{[0,t_1]}}]{1}{}} 
     ~=    -\sum_l\norm{f}{\Lp[\GammaT]{\infty}{}}\norm{\sqrt{\cl+a}}{\Lp[{\Gamma\times{[0,t_1]}}]{2}{}}^2 \\
     &\geq -\delta\sum_l\norm{\grad\sqrt{\cl+a}}{\Lp[{\Omega\times{[0,t_1]}}]{2}{}}^2 -2\delta^{-1}\norm{f}{\Lp[\GammaT]{\infty}{}}^2\sum_l\norm{\sqrt{\cl+a}}{\Lp[{\Omega\times{[0,t_1]}}]{2}{}}^2 \\
     &=    -\delta\sum_l\norm{\grad\sqrt{\cl+a}}{\Lp[{\Omega\times{[0,t_1]}}]{2}{}}^2 -2\delta^{-1}\norm{f}{\Lp[\GammaT]{\infty}{}}^2\sum_l\norm{\cl+a}{\Lp[{\Omega\times{[0,t_1]}}]{1}{}} \\
     &\geq -\delta\sum_l\norm{\grad\sqrt{\cl+a}}{\Lp[{\Omega\times{[0,t_1]}}]{2}{}}^2 -2\delta^{-1}\norm{f}{\Lp[\GammaT]{\infty}{}}^2\sum_l\norm{\lyapunov(\cl+a)}{\Lp[{\Omega\times{[0,t_1]}}]{1}{}} . 
\end{align*}
This yields for the sum of the convection integrals
\begin{align*}
 I_{co} &\geq -\delta\sum_l\norm{\grad\sqrt{\cl+a}}{\Lp[{\Omega\times{[0,t_1]}}]{2}{}}^2 -2\delta^{-1}\norm{f}{\Lp[\GammaT]{\infty}{}}^2\sum_l\norm{\lyapunov(\cl+a)}{\Lp[{\Omega\times{[0,t_1]}}]{1}{}} \\
        &~~~~ -\sum_l\scp{f}{a\ln(\cl+a)}_{\Gamma\times{[0,t_1]}}.
\end{align*}
Similarly, we transform the sum of the electric drift integrals with integration by parts and equation~\eqref{eq:gaussWeakDIV-fpi} to
\begin{align*}
 I_{el}:=&     -\sum_l\scp{\cl\coeffEL\fieldEL}{\grad\ln(\cl+a)}_{\Omega\times{[0,t_1]}}  \\
         &=    -e(\permitELscalar\boltz\temp)^{-1}\sum_l\scp{\fieldEL}{\zl\grad\cl}_{\Omega\times{[0,t_1]}}  - a\sum_l\coeffEL\scp{\fieldEL}{\zl\grad\ln(\cl+a)}_{\Omega\times{[0,t_1]}}  \\[1.0mm]
         &=    e(\permitELscalar\boltz\temp)^{-1}\scp{\chargeEL}{\chargeEL}_{\Omega\times{[0,t_1]}}   ~-~  e(\permitELscalar\boltz\temp)^{-1}\scp{\sigma}{\chargeEL}_{\Gamma\times{[0,t_1]}} \\
         &~~~~ +\sum_l\coeffEL\sqbrac{~\scp{\chargeEL}{a\ln(\cl+a)}_{\Omega\times{[0,t_1]}} -\scp{\sigma}{a\ln(\cl+a)}_{\Gamma\times{[0,t_1]}}~} \\
         &=:~B.1 +B.2 +B.3.
\end{align*}
Exactly as we treated the integral~$A.1$, we come for the integral~$B.2$ to
\begin{align*}
 B.2 &=     -e(\permitELscalar\boltz\temp)^{-1} \sum_l \scp{\sigma}{\zl\cl}_{\Gamma\times{[0,t_1]}} 
     ~\geq  -\frac{e\max_l\abs{\zl}}{\permitELscalar\boltz\temp} \sum_l \scp{\abs{\sigma}}{\abs{\cl+a} }_{\Gamma\times{[0,t_1]}} \\
     &\geq  -\frac{e\max_l\abs{\zl}}{\permitELscalar\boltz\temp}\norm{\sigma}{\Lp[\GammaT]{\infty}{}}\norm{\sqrt{\cl+a}}{\Lp[{\Gamma\times{[0,t_1]}}]{2}{}}^2\\
     &\geq -\delta\sum_l\norm{\grad\sqrt{\cl+a}}{\Lp[{\Omega\times{[0,t_1]}}]{2}{}}^2 -\frac{2e^2\max_l\abs{\zl}^2}{\delta(\permitELscalar\boltz\temp)^2}\norm{\sigma}{\Lp[\GammaT]{\infty}{}}^2\sum_l\norm{\lyapunov(\cl+a)}{\Lp[{\Omega\times{[0,t_1]}}]{1}{}} .
\end{align*}
Hence, we arrive for the sum of the electric drift integrals at
\begin{align*}
I_{el} &\geq e(\permitELscalar\boltz\temp)^{-1} \norm{\chargeEL}{\Lp[{\Omega\times{[0,t_1]}}]{2}{}}^2 
             -\frac{2e^2\max_l\abs{\zl}^2}{\delta(\permitELscalar\boltz\temp)^2}\norm{\sigma}{\Lp[\GammaT]{\infty}{}}^2\sum_l\norm{\lyapunov(\cl+a)}{\Lp[\Omega\times{[0,t_1]}]{1}{}} \\
       &~~   -\delta\sum_l\norm{\grad\sqrt{\cl+a}}{\Lp[\Omega\times{[0,t_1]}]{2}{}}^2 
             + \sum_l\frac{e\zl}{\permitELscalar\boltz\temp}\sqbrac{~\scp{\chargeEL}{a\ln(\cl+a)}_\Omega-\scp{\sigma}{a\ln(\cl+a)}_\Gamma~}.
\end{align*}
The remaining sum of the reaction integrals, we estimate with \cref{lemma:lyapunov} and $\cl\geq0$ according to \cref{lemma:nonnegative-Lcomp} by
\begin{align*}
 &     \sum_l \scp{\Rl(\Cl)}{\ln(\cl+a)}_{\Omega\times{[0,t_1]}}
 ~\leq \max_l C_{\Rl} \sum_l \scp{\frac{\abs{\Cl}}{(\cl+a)}}{(\cl+a)\ln(\cl+a)}_{\Omega\times{[0,t_1]}} \\
 &\leq \max_l C_{\Rl} \sum_l \scp{\frac{\abs{\Cl}}{(\cl+a)}}{\lyapunov(\cl+a)}_{\Omega\times{[0,t_1]}} \\
 &\leq \max_l C_{\Rl} \scp{\frac{\abs{\Cl}}{\sum_l(\cl+a)}}{\sum_k\lyapunov(\cl[k]+a)}_{\Omega\times{[0,t_1]}}  
 ~\leq    \max_l C_{\Rl} \sum_l\norm{\lyapunov(\cl+a)}{\Lp[\Omega\times{[0,t_1]}]{1}{}} .
\end{align*}
Combining the preceding estimates leads with a suitable choice of the free parameter~$\delta>0$ 
to the intermediate inequality 
\begin{align*}
 &     \sum_l \norm{\lyapunov(\cl(t_1)+a)}{\Lp{1}{}} + \frac{\alpha_D}{2} \sum_l \norm{\grad\sqrt{\cl+a}}{\Lp[{\Omega\times{[0,t_1]}}]{2}{}}^2 +\frac{e}{\permitELscalar\boltz\temp} \norm{\chargeEL}{\Lp[{\Omega\times{[0,t_1]}}]{2}{}}^2 \\
 &\leq \sum_l \norm{\lyapunov(\clstart+a)}{\Lp{1}{}} + \frac{8}{\alpha_D}\norm{f}{\Lp[\GammaT]{\infty}{}}^2\sum_l\norm{\lyapunov(\cl+a)}{\Lp[{\Omega\times{[0,t_1]}}]{1}{}} \\
 &~~~~ +\frac{8e^2\max_l\abs{\zl}^2}{\alpha_D(\permitELscalar\boltz\temp)^2}\norm{\sigma}{\Lp[\GammaT]{\infty}{}}^2\sum_l\norm{\lyapunov(\cl+a)}{\Lp[\Omega\times{[0,t_1]}]{1}{}} + \sum_l\scp{f}{a\ln(\cl+a)}_{\Gamma\times{[0,t_1]}} \\
 &~~~~ -\sum_l\coeffEL\sqbrac{~\scp{\chargeEL}{a\ln(\cl+a)}_{\Omega\times{[0,t_1]}} -\scp{\sigma}{a\ln(\cl+a)}_{\Gamma\times{[0,t_1]}}~} \\
 &~~~~ + \max_l C_{\Rl} \sum_l\norm{\lyapunov(\cl+a)}{\Lp[\Omega\times{[0,t_1]}]{1}{}} .
\end{align*}
We now can safely let $a\searrow0$. Thereby, we obtain with $a\ln(\cl+a)\rightarrow0$ the entropy estimate
\begin{align*}
 &     \sum_l \norm{\lyapunov(\cl(t_1))}{\Lp{1}{}} + \frac{\alpha_D}{2} \sum_l \norm{\grad\sqrt{\cl}}{\Lp[{\Omega\times{[0,t_1]}}]{2}{}}^2 +\frac{e}{\permitELscalar\boltz\temp} \norm{\chargeEL}{\Lp[{\Omega\times{[0,t_1]}}]{2}{}}^2 \\
 &\leq \sum_l \norm{\lyapunov(\clstart)}{\Lp{1}{}} \\
 &~~~~ +\underbrace{\frac{8}{\alpha_D}\norm{f}{\Lp[\GammaT]{\infty}{}}^2 +\frac{8e^2\max_l\abs{\zl}^2}{\alpha_D(\permitELscalar\boltz\temp)^2}\norm{\sigma}{\Lp[\GammaT]{\infty}{}}^2 + \max_l C_{\Rl} }_{=:b} \sum_l\norm{\lyapunov(\cl)}{\Lp[{\Omega\times{[0,t_1]}}]{1}{}} .
\end{align*}
Applying \Grownwall\ leads immediately to
\begin{align*}
 \sum_l\norm{\lyapunov(\cl)}{\Lp[I]{\infty}{;\Lp{1}{}}} \leq \sqbrac{1+ b\timeEnd e^{b\timeEnd}}\sum_l \norm{\lyapunov(\clstart)}{\Lp{1}{}}=:\tilde{C}_L(\timeEnd) .
\end{align*}
Hence, we obtain
\begin{align*}
 \norm{\chargeEL}{\Lp[{\Omega\times{[0,t_1]}}]{2}{}}^2 \leq  \frac{\permitELscalar\boltz\temp}{e} \sqbrac{1+ b\tilde{C}_L(\timeEnd)} \sum_l\norm{\lyapunov(\clstart)}{\Lp{1}{}} =: C_{L,1}.
\end{align*}
From this bound, we deduce with \cref{lemma:lyapunov}
\begin{align*}
 &      \sum_l\norm{\sqrt{\cl}}{\Lp[I]{\infty}{;\Lp{2}{}}}^2 + \frac{\alpha_D}{2} \sum_l \norm{\grad\sqrt{\cl}}{\Lp[{\Omega\times{[0,t_1]}}]{2}{}}^2 \\
 &\leq  \frac{2}{\alpha_D} \sqbrac{1+ b\tilde{C}_L(\timeEnd)} \sum_l\norm{\lyapunov(\clstart)}{\Lp{1}{}} =: C_{L,2}.
\end{align*}
Finally, we define the constant $C_L=C_L(\timeEnd)$ by $C_L:=2\max(C_{L,1},C_{L,2})$.\\[2.0mm]
%
\underline{Case 2: $\cl\neq\clold$}~~
Concerning the first stated estimate, we assume that $\clold$ is contained in a ball in $\Lp[\OmegaT]{2}{}$ with radius $R$. Thus, we trivially obtain 
\begin{align*}
  \norm{\chargeELold}{\Lp[\OmegaT]{2}{}}^2\leq \max_l\abs{\zl}\sum_l \norm{\clold}{\Lp[\OmegaT]{2}{}}^2  \leq LR^2\max_l\abs{\zl}~.
\end{align*}
Concerning the second stated estimate, we note that the bound of the time integrals, the diffusion integrals, the convection integrals, and the reaction integrals remain the same as in case~1 above. 
Only in the estimate for electric drift integral, the bound for the subintegral~$B.1$ changes. More precisely, this time we get for the integral~$B.1$ with \ref{lemma:lyapunov}
\begin{align*}
 B.1 &\geq  -e(\permitELscalar\boltz\temp)^{-1}\max_l\abs{\zl} \norm{\chargeELold}{\Lp[\OmegaT]{\infty}{}} \sum_l\norm{\cl}{\Lp[\Omega\times{[0,t_1]}]{1}{}} \\
     &\geq  -e(\permitELscalar\boltz\temp)^{-1}\max_l\abs{\zl} \norm{\chargeELold}{\Lp[\OmegaT]{\infty}{}} \sum_l\norm{\lyapunov(\cl+a)}{\Lp[\Omega\times{[0,t_1]}]{1}{}}~. 
\end{align*}
Analogously to case~1, we obtain after $a\searrow0$ and $a\ln(\cl+a)\rightarrow0$ the entropy estimate
\begin{align*}
 &     \sum_l \norm{\lyapunov(\cl(t_1))}{\Lp{1}{}} + \frac{\alpha_D}{2} \sum_l \norm{\grad\sqrt{\cl}}{\Lp[{\Omega\times{[0,t_1]}}]{2}{}}^2 \\
 &\leq \sum_l \norm{\lyapunov(\clstart)}{\Lp{1}{}} +\underbrace{e(\permitELscalar\boltz\temp)^{-1}\max_l\abs{\zl} \norm{\chargeELold}{\Lp[\OmegaT]{\infty}{}} }_{=:b_0} \sum_l\norm{\lyapunov(\cl)}{\Lp[{\Omega\times{[0,t_1]}}]{1}{}}\\
 &~~~~ +\underbrace{\frac{8}{\alpha_D}\norm{f}{\Lp[\GammaT]{\infty}{}}^2 +\frac{8e^2\max_l\abs{\zl}^2}{\alpha_D(\permitELscalar\boltz\temp)^2}\norm{\sigma}{\Lp[\GammaT]{\infty}{}}^2 + \max_l C_{\Rl} }_{=:b} \sum_l\norm{\lyapunov(\cl)}{\Lp[{\Omega\times{[0,t_1]}}]{1}{}} .
\end{align*}
Provided that $\clold$ is contained in a ball in $\Lp[\OmegaT]{\infty}{}$ with radius $R$, we have with \Grownwall\ 
\begin{align*}
 \sum_l\norm{\lyapunov(\cl)}{\Lp[I]{\infty}{;\Lp{1}{}}} \leq \sqbrac{1+ (b+b_0)\timeEnd e^{(b+b_0)\timeEnd}}\sum_l \norm{\lyapunov(\clstart)}{\Lp{1}{}}=:\tilde{C}_L(\timeEnd,R) .
\end{align*}
Thereby, we finally arrive at 
\begin{align*}
 &      \sum_l\norm{\sqrt{\cl}}{\Lp[I]{\infty}{;\Lp{2}{}}}^2 + \frac{\alpha_D}{2} \sum_l \norm{\grad\sqrt{\cl}}{\Lp[{\Omega\times{[0,t_1]}}]{2}{}}^2 \\
 &\leq  \frac{2}{\alpha_D} \sqbrac{1+ (b+b_0)\tilde{C}_L(\timeEnd,R)} \sum_l\norm{\lyapunov(\clstart)}{\Lp{1}{}} =: C_L(\timeEnd,R).
\end{align*}
\end{proof}
Next, we prove the crucial energy estimates for the chemical species~$\cl$ with the aid of the entropy estimates from \cref{lemma:entropyBound}.  
%
%
\begin{lemma}[Energy estimates]\label{lemma:energy-Lcomp}
Let \ref{Assump:Geom}--\ref{Assump:back-charge} be valid and let $\brac{\fieldEL,\potEL,\fieldF,\press,\Cl}\in \setR^{2+2n+L}$ be a weak solution 
of the system~\eqref{eq:Model-1a}--\eqref{eq:Model-1j} according to \cref{def:FixOp}. Then, we have 
\begin{align*}
  \sum_l\sqbrac{ \norm{\cl}{\Lp[I]{\infty}{;\Lp{2}{}}}+\norm{\cl}{\Lp[I]{2}{;\spaceTC}} }~\leq~ C_0
\end{align*}
Herein, the dependency of the constant is
\begin{align*}
 C_0 &= C_0\!\brac{C_L,\timeEnd,\max_l\abs{\zl},\norm{f}{\Lp[\GammaT]{\infty}{}},\norm{\sigma}{\Lp[\GammaT]{\infty}{}},\norm{\rho_b}{\Lp[\OmegaT]{\infty}{}},\norm{\clstart}{\Lp{2}{}} }.
\end{align*}
\end{lemma}
\begin{proof}
For ease of readability, we split the proof into two cases. Moreover, in the current proof we use \GagNirenberg\ in way, which restricts this result to two space dimensions.\\[2.0mm]
\underline{Case 1: $n=2$ and $\cl=\clold$}~~
In equations~\eqref{eq:transportWeak-fpi}, we choose the test functions $\varphi := \cl \in \spaceTC$ and we sum over $l=1,2$. Thereby, we get for the time integrals and 
the diffusion integrals with \ref{Assump:Ellip}
\begin{align*}
      \sum_l \dualp{\theta\dert\cl}{\cl}_{1,\Omega} + \sum_l\scp{\Dl\grad\cl}{\grad\cl}_\Omega 
 \geq \frac{\theta}{2}\derr\sum_l \norm{\cl}{\Lp{2}{}}^2 + \alpha_D \sum_l\norm{\grad\cl}{\Lp{2}{}}^2~~.
\end{align*}
For the convection integrals, we firstly use integration by parts and we insert equation \eqref{eq:darcyWeakDIV-fpi}. Secondly, we use \Holder\ and \cref{lemma:interpol-boundary} 
with a rescaled free parameter $\delta=\norm{f}{\Lp[\GammaT]{\infty}{}}^{-1}\delta$. This leads us to
\begin{align*}
 &     -\sum_l\scp{{\cl}\fieldF}{\grad\cl}_\Omega 
  =    -\frac{1}{2} \sum_l \scp{\fieldF}{\grad\cl^2}_\Omega 
  =    -\frac{1}{2} \sum_l\scp{f}{\cl^2}_\Gamma 
  \geq -\frac{1}{2} \sum_l \norm{f}{\Lp[\GammaT]{\infty}{}} \norm{\cl}{\Lp[\Gamma]{2}{}}^2 \\
 &\geq -\delta \sum_l \norm{\grad\cl}{\Lp{2}{}}^2 
  ~-~  \delta^{-1} \norm{f}{\Lp[\GammaT]{\infty}{}}^2 \sum_l \norm{\cl}{\Lp{2}{}}^2 .
\end{align*}
Analogously, for the electric~drift integral we firstly integrate by parts, we insert equation \eqref{eq:gaussWeakDIV-fpi}. Then, we use \Holder\ and \cref{lemma:interpol-boundary}. This yields
\begin{align*}
 I_{el}&:=   -\frac{e}{\permitELscalar\boltz\temp} \sum_l\scp{\zl{\cl}\fieldEL}{\grad\cl}_\Omega
       ~=    \frac{e}{\permitELscalar\boltz\temp} \sum_l\frac{\zl}{2}\sqbrac{~\scp{\rho_b+\chargeELold}{\cl^2}_\Omega - \scp{\sigma}{\cl^2}_\Gamma~}\\
       &\geq \frac{e}{\permitELscalar\boltz\temp} \sum_l\frac{\zl}{2}\scp{\chargeEL}{\cl^2}_\Omega 
             -\delta \sum_l\abs{\zl}\norm{\grad\cl}{\Lp{2}{}}^2 \\
       &~~~~ -\frac{2e^2}{(\permitELscalar\boltz\temp)^2}\max_l\abs{\zl} \sqbrac{\delta^{-1}\norm{\sigma}{\Lp[\GammaT]{\infty}{}}^2+\norm{\rho_b}{\Lp[\OmegaT]{\infty}{}}} \sum_l\abs{\zl}\norm{\cl}{\Lp{2}{}}^2 \\
       &~~~~ := \labelInt + I.b + I.c. 
\end{align*}
We now have to bound the integral~$\labelInt$. For that purpose, we apply \Holder, \GagNirenberg\ (we have $n=2$) and \Young, which yields
\label{page:signCondition}
\begin{align*}
 \labelInt &\leq \frac{e}{\permitELscalar\boltz\temp} \max_l\abs{\zl} \norm{\chargeEL}{\Lp{2}{}} \sum_l\norm{\cl^2}{\Lp{2}{}} 
            =    \frac{e}{\permitELscalar\boltz\temp} \max_l\abs{\zl} \norm{\chargeEL}{\Lp{2}{}} \sum_l\norm{\cl}{\Lp{4}{}}^2 \\
           &\leq \frac{e}{\permitELscalar\boltz\temp} \max_l\abs{\zl} \norm{\chargeEL}{\Lp{2}{}} \sum_l\norm{\cl}{\Lp{2}{}}\norm{\cl}{\Hk{1}{}} \\
           &\leq \delta \sum_l\norm{\cl}{\Hk{1}{}}^2 +  \frac{e^2}{\delta (\permitELscalar\boltz\temp)^2} \max_l\abs{\zl}^2 \norm{\chargeEL}{\Lp{2}{}}^2 \sum_l\norm{\cl}{\Lp{2}{}}^2~.
\end{align*}
Substituting the bound for $\labelInt$, leads for electric drift integral to
\begin{flalign*}
&~~I_{el} \geq -\delta \sum_l\norm{\grad\cl}{\Lp{2}{}}^2 
             -\frac{2e^2\max_l\abs{\zl}^2}{\delta(\permitELscalar\boltz\temp)^2} \sqbrac{\norm{\sigma}{\Lp[\GammaT]{\infty}{}}^2+\norm{\rho_b}{\Lp[\OmegaT]{\infty}{}}+\norm{\chargeEL}{\Lp{2}{}}^2} \sum_l\norm{\cl}{\Lp{2}{}}^2 . &
\end{flalign*}
For the reaction integrals, we applying \ref{Assump:Reaction}, \Young, and we recall $\cl\geq0$. This results for the reaction integrals in
\begin{align*}
  &    \sum_l\scp{\theta\Rl(\Cl)}{\cl}_\Omega 
   \leq \theta \max_l C_{\Rl} \scp{\abs{\Cl}}{\sum_l\cl}_\Omega 
   \leq \theta \max_l C_{\Rl} \scp{\sum_l\cl}{\sum_l\cl}_\Omega \\
  &=    \theta \max_l C_{\Rl} \sum_l \scp{\cl}{\cl}_\Omega
        +\theta \max_l C_{\Rl} \sum_{k\neq l} \scp{\cl[k]}{\cl}_\Omega 
  \leq 3\theta C_{\Rl} \sum_l \norm{\cl}{\Lp{2}{}}^2 ~.
\end{align*}
By combining the preceding estimates, we deduce with the choice~$\delta:=\frac{\alpha_D}{4}$ the estimate
\begin{align*}
 &    \frac{\theta}{2}\derr\sum_l \norm{\cl}{\Lp{2}{}}^2 + \frac{\alpha_D}{2} \sum_l\norm{\grad\cl}{\Lp{2}{}}^2 \nonumber\\
 &\leq  \frac{12e^2\max_l\abs{\zl}^2}{\alpha_D(\permitELscalar\boltz\temp)^2}\sqbrac{\norm{\sigma}{\Lp[\GammaT]{\infty}{}}^2+\norm{\rho_b}{\Lp[\OmegaT]{\infty}{}}+ \norm{\chargeEL}{\Lp{2}{}}^2} \sum_l\norm{\cl}{\Lp{2}{}}^2  \nonumber\\
 &~~~~  +\sqbrac{\frac{4}{\alpha_D} \norm{f}{\Lp[\GammaT]{\infty}{}}^2 +  3\theta \max_lC_{\Rl} }\sum_l\norm{\cl}{\Lp{2}{}}^2~. 
\end{align*}
With the abbreviations   
\begin{align*}
 b(t) &:=  \frac{12e^2\max_l\abs{\zl}^2}{\alpha_D(\permitELscalar\boltz\temp)^2}\norm{\chargeEL}{\Lp{2}{}}^2 \\
  A_0 &:=    \min\!\brac{\frac{\theta}{2},\frac{\alpha_D}{2}} \\
 B_0  &:=  \frac{12e^2\max_l\abs{\zl}^2}{A_0\alpha_D(\permitELscalar\boltz\temp)^2} \left[ \norm{\sigma}{\Lp[\GammaT]{\infty}{}}^2
                                                                 +\norm{\rho_b}{\Lp[\OmegaT]{\infty}{}} 
      + \norm{f}{\Lp[\GammaT]{\infty}{}}^2 +  3\theta \max_lC_{\Rl}  \right]~,
\end{align*}
we deduce from the preceding estimate  with \Grownwall\ and \cref{lemma:entropyBound}
\begin{align*}
       \sum _l\norm{\cl}{\Lp[I]{\infty}{;\Lp{2}{}}}^2 
 &\leq \exp\!\brac{ \Intdt[][\timeEnd]{B_0 +b(t)~} } \sum_l\norm{\clstart}{\Lp{2}{}}^2 \\
 &\leq e^{B_0\timeEnd} \exp\!\brac{ \frac{12e^2\max_l\abs{\zl}^2}{\alpha_D(\permitELscalar\boltz\temp)^2} \norm{\chargeEL}{\Lp[\OmegaT]{2}{}}^2 } \sum_l\norm{\clstart}{\Lp{2}{}}^2 \\
 &\leq e^{B_0\timeEnd} \exp\!\brac{ \frac{12e^2\max_l\abs{\zl}^2}{\alpha_D(\permitELscalar\boltz\temp)^2} C_L^2 } \sum_l\norm{\clstart}{\Lp{2}{}}^2
 ~:=   \hat{C}^2_0(\timeEnd).
\end{align*}
Substituting this bound into the above estimate   and integrating in time over $[0,\timeEnd]$, yields the desired a~priori estimate 
\begin{align*}
        \sum _l \sqbrac{\norm{\cl}{\Lp[I]{\infty}{;\Lp{2}{}}} + \norm{\grad\cl}{\Lp[\OmegaT]{2}{}} } 
  &\leq \hat{C}_0 + B_0^{\frac{1}{2}}\timeEnd^{\frac{1}{2}}\sum_l\norm{\cl}{\Lp[I]{\infty}{;\Lp{2}{}}} \nonumber\\
  &\leq \hat{C}_0 + B_0^{\frac{1}{2}}\timeEnd^{\frac{1}{2}}\hat{C}_0 
   ~:=  C_0(\timeEnd)~.
\end{align*}
%
%
\underline{Case 2: $n=2$ and $\cl\neq\clold$}~~
Again, we test equations~\eqref{eq:transportWeak-fpi} with $\test:=\cl\in\spaceTC$ and we run through the same calculations as already carried out in Case~1. The only difference is 
that we now immediately obtain for the integral~$\labelInt$ 
\begin{align*}
 \labelInt  \leq \frac{e}{\permitELscalar\boltz\temp} \max_l\abs{\zl} \norm{\chargeELold}{\Lp{\infty}{}} ~\sum_l\norm{\cl}{\Lp{2}{}}^2 
            \leq \frac{e}{\permitELscalar\boltz\temp} \max_l\abs{\zl}^2 \norm{\Clold}{\Lp{\infty}{}}~~ \sum_l\norm{\cl}{\Lp{2}{}}^2
\end{align*}
We note, that in \cref{def:FixOp} we introduced the space~$X$ and the set~$\fixSet\subset X$. Furthermore, we supposed $\Clold\in\fixSet$, which ensures that the $L^\infty$-norms of the $\clold$ remain finite. 
Thus, provided we know $\norm{\Clold}{\Lp{\infty}{}}\leq R$ for all $\Clold\in\fixSet$, the constant in the above estimate just depends on an additional parameter~$R$. 
In conclusion, we analogously obtain 
\begin{align*}
       \sum _l \sqbrac{\norm{\cl}{\Lp[I]{\infty}{;\Lp{2}{}}} + \norm{\grad\cl}{\Lp[\OmegaT]{2}{}} } 
  \leq \hat{C}_0 + B_0^{\frac{1}{2}}\timeEnd^{\frac{1}{2}}\hat{C}_0 
   ~:=  C_0(\timeEnd,\norm{\Clold}{\Lp{\infty}{}})~.
\end{align*}
Herein, we just changed the definition of the constant~$b(t)$ by
\begin{align*}
 b(t) :=  \frac{12e^2\max_l\abs{\zl}^2}{\alpha_D(\permitELscalar\boltz\temp)^2}\norm{\Clold}{\Lp{\infty}{}} 
      \leq \frac{12e^2\max_l\abs{\zl}^2}{\alpha_D(\permitELscalar\boltz\temp)^2}\norm{\Clold}{\Lp[\OmegaT]{\infty}{}} =:b_0~.
\end{align*}
\end{proof}
%
%
\subsection{Existence of a fixed point}\label{subsec:fpi-point}
For the remaining part of the proof, we can shorty refer to Part~1 of this work, as we wrote Part~1 such that we can direct copy the following parts of the proof. 
The only difference is, that we subsequently consider $L$ solutes instead of two solutes. Thus, the range of the index~$l$ is now $l=1,\ldots,L$ instead of $l=1,2$. 
%
%
\begin{lemma}[Boundedness]\label{lemma:bounded-Lcomp}
Let \ref{Assump:Geom}--\ref{Assump:back-charge} be valid and let $\brac{\fieldEL,\potEL,\fieldF,\press,\Cl}\in \setR^{2+2n+L}$ be a weak solution of the system~\eqref{eq:Model-1a}--\eqref{eq:Model-1j} 
according to \cref{def:FixOp}. Then, we have  
\begin{align*}
  \sum_l \norm{\cl}{\Lp[\OmegaT]{\infty}{}} ~\leq~ C_M~.
\end{align*}
Herein, the dependency of the constant is
\begin{align*}
 C_M =C_M \!\brac{C_0, \timeEnd,\max_l\abs{\zl}, \norm{f}{\Lp[\GammaT]{\infty}{}},\norm{\sigma}{\Lp[\GammaT]{\infty}{}},\norm{\rho_b}{\Lp[\OmegaT]{\infty}{}},\norm{\clstart}{\Lp{\infty}{}} }.
\end{align*}
\end{lemma}
\begin{proof}
The proof is identical to the proof of the corresponding proof of Part~1 of this work.
\end{proof} 
%
%
\begin{thm}[A priori Bounds]\label{thm:aprioriBounds-Lcomp}
Let \ref{Assump:Geom}--\ref{Assump:back-charge} be valid and let $\brac{\fieldEL,\potEL,\fieldF,\press,\Cl}\in \setR^{2+2n+L}$ be a weak solution 
of \eqref{eq:Model-1a}--\eqref{eq:Model-1j} according to \cref{def:weaksolution}. Then, we have
\begin{align*}
 & \norm{\potEL}{\Lp[I]{\infty}{;\Lp{2}{}}} + \norm{\fieldEL}{\Lp[I]{\infty}{;\Lp{2}{}}} ~\leq~ C(\timeEnd)~, \\
 & \norm{\press}{\Lp[I]{\infty}{;\Lp{2}{}}} + \norm{\fieldF}{\Lp[I]{\infty}{;\Lp{2}{}}}  ~\leq~ C(\timeEnd)~, \\
 & \sum_l\sqbrac{ \norm{\cl}{\Lp[I]{\infty}{;\Lp{2}{}}}+\norm{\cl}{\Lp[I]{2}{;\spaceTC}}+\norm{\cl}{\Hk[I]{1}{;\spaceTC^\ast}}+\norm{\cl}{\Lp[\OmegaT]{\infty}{}} }~\leq~ C(\timeEnd)~.
\end{align*}
\end{thm}
\begin{proof} 
The proof is identical to the corresponding proof of Part~1 of this work. 
\end{proof}
%
%
\begin{thm}\label{thm:existence-Lcomp}
Let \ref{Assump:Geom}--\ref{Assump:back-charge} be valid. Then, there exists a solution $\brac{\fieldEL,\potEL,\fieldF,\press,\Cl}\in \setR^{2+2n+L}$ of 
equations~\eqref{eq:Model-1a}--\eqref{eq:Model-1j} according to \cref{def:weaksolution}.
\end{thm}
\begin{proof}
The proof is identical to the proof of the corresponding proof of Part~1 of this work.
\end{proof}
%
%
\section{Conclusion}
In this paper, we showed the global existence of unique solutions of the \dpnp. The contribution of this paper was to deal with multicomponent electrolyte solutions, 
which consist of a neutral solvent and multiple charged solutes~$\cl$ with arbitrary valencies. In this situation, the main difficulty was to establish a~priori estimates for 
the chemical species~$\cl$. By using entropy estimates, which have been developed for certain Lyapunov functionals, we successfully obtained such a~priori estimates. 
However, these techniques are restricted to two space dimensions, as we combined the entropy estimates with \GagNirenberg\ in a way, which is valid only in two space dimensions. 
\par
In particular, by means of the entropy estimates the presented proof avoids further restrictions on the electrolyte solutions, such as the often used electroneutrality constraint or the volume additivity constraint.
Therefore, our results can be applied to two dimensional models of general electrolyte solutions, which are captured by the \dpnp. This important especially in biological applications and hydrodynamical applications.   
\par
Finally, we note that if we restrict us to electrolyte solutions at rest (in this case the equations~\eqref{eq:Model-1e}--\eqref{eq:Model-1g} coming from Darcy's law vanish) 
and linear reaction rates, the presented model in this paper is identical to the considered model in \cite{Glitzky2004}, if additionally no constraints are involved in \cite{Glitzky2004}. 
However, even in this situation we gave a new proof of the same result. More precisely, in \cite{Glitzky2004} the crucial a~priori estimates were obtained by involving a two dimensional, nonlinear version 
of \GagNirenberg\ from \cite{biler1994}, whereas we used the entropy estimates for that purpose.
%
%
%
\section*{Acknowledgements}
M. Herz is supported by the Elite Network of Bavaria.  
%
%
%
\printbibliography
\end{document}